\documentclass[11pt]{llncs}
\usepackage{graphics}
\def\Q{{\bf Q}}
\begin{document}
\title{Elliptic curves of large rank and small conductor}
\author{Noam D. Elkies\inst{1}\thanks{Supported in part
  by NSF grant DMS-0200687.}
\and Mark Watkins\inst{2}\thanks{Supported in part
  by VIGRE postdoctoral funding from the NSF,
  and was a Visiting Scholar with the MAGMA Computer Algebra Group
  at the University of Sydney
  during part of the time in which this work was done.}}
\institute{Department of Mathematics, Harvard University,
Cambridge, MA 02138\\
\email{elkies@math.harvard.edu}\\
\and
The Pennsylvania State University\\
\email{watkins@math.psu.edu}\\
}
\maketitle
\begin{abstract}
For $r=6,7,\ldots,11$ we find an elliptic curve $E/\Q$
of rank at least~$r$ and the smallest conductor known,
improving on the previous records by factors ranging from
$1.0136$ (for $r=6$) to over $100$ (for $r=10$ and $r=11$).
We describe our search methods, and tabulate,
for each $r=5,6,\ldots,11$, the five curves of lowest conductor,
and (except for $r=11$) also the five of lowest absolute discriminant,
that we found.
\end{abstract}
\section{Introduction and Motivation}
An elliptic curve over the rationals is a curve~$E$\/ of genus~$1$,
defined over~$\Q$, together with a \hbox{$\Q$-rational} point.
A theorem of Mordell \cite{mordell} states that the
rational points on~$E$\/ form a finitely generated abelian group under a
natural group law. The rank of~$E$\/ is the rank of the free
part of this group. Currently there is no general unconditional algorithm
to compute the rank. Elliptic curves of large rank are hard to find;
the current record is a curve of rank at least~$24$
(see \cite{rank24}).\footnote{Conrey wrote of a curve of rank~$26$
\cite[p.~353]{rank26}, but confirms in e-mail to the authors that ``$26$''
was a typographical error for ``$24$'' as in~\cite{rank24}.}

We investigate a slightly different question: instead of seeking
curves of large rank, we fix a small rank~$r$ (here $5 \leq r \leq 11$)
and try to make the conductor~$N$\/ as small as possible,
which, due to the functional equation for the $L$-function
of the elliptic curve, is more natural than trying to minimize
the absolute discriminant~$|\Delta|$.
The question of how fast the rank can grow as a function of~$N$\/
has generated renewed interest lately, partially due to the predictions made
by random matrix theory about $\zeta$-function analogues~\cite{conrey-gonek}.
However, there are at present two different conjectures,
one that comes from a function field analogue,
and another from analytic number theory considerations.
We shall return to this in Section~\ref{rankcond}.

We try to find $E/\Q$\/ of high rank and low
conductor by searching for elliptic curves that have many integral points.
As stated, this strategy is ill-posed, as integrality of points is not
invariant under change of model (defining equation).
However, if we only consider (say) N\'eron models
then the question makes sense,
and a conjecture of Lang \cite[p.~140]{lang}
links the number of integral points to the rank,
at least in one direction.
More explicitly, one might conjecture that there is an
absolute constant~$C$\/ such that the number of integral points
on an elliptic curve~$E$\/ of rank~$r$ is bounded by~$C^{r+1}$.
The best result is due to Silverman~\cite{silverman},
who shows the conjecture is true when the $j$-invariant $j(E)$ is integral,
and in fact proves that for every number field~$K$\/ there is a
constant~$C_K$ such that the number of \hbox{$S$-integral} points over~$K$\/
is bounded by $C_K^{(1+r)(1+\delta)+|S|}$, where $\delta$
is the number of primes of $K$\/ at which $j(E)$ is nonintegral.
Explicit constants appear in~\cite{gross-silverman}.
Szpiro's conjecture~\cite{szpiro}, which is equivalent to the
Masser-Oesterl\'e ABC conjecture~\cite{ABC},
states that~$\Delta \ll N^{6+o(1)}$;
Hindry and Silverman \cite{hindry-silverman} show this implies
that the number of \hbox{$S$-integral} points on a quasi-minimal model
of~$E/K$ is bounded by $C_K^{(1+r)\sigma_{E/K}+|S|}$ where $\sigma_{E/K}$ is
the Szpiro ratio, which is the ratio of the logarithms of the norms
of the discriminant and the conductor of~$E/K$\/.
Finally, Abramovich~\cite{abramovich} has shown that
the Lang-Vojta conjecture (which states that the integral points
on a quasi-projective variety of log general type are not Zariski dense,
see~\cite[4.4]{vojta}) implies the uniform boundedness of the number of
integral points on rational semistable elliptic curves,
but the lack of control over the Zariski closure of the integral points
makes this result ineffective.

Conversely, it is frequently the case that elliptic curves of high rank,
and especially those with relatively small conductor, have many
integral points, and thus our search method is likely to find
these curves.  In fact, for each $r$ in our range $5 \leq r \leq 11$
we found a curve~$E$\/ of rank at least~$r$ whose conductor~$N$\/
is the smallest known.  For $r=5$ this was a previously known
(see~\cite{brumer-mcguinness}) curve with $N=19047851$.
For the other~$r$ our curve is new,
with $N$\/ smaller than the previous record by a factor ranging from
$1.0136$ for $r=6$ to over $100$ for $r=10$ and $r=11$.
As a byproduct we also find the curves of rank~$r$
whose discriminants $\Delta$ have the smallest absolute values known.
We estimate that finding a similarly good rank~12 curve would take
\hbox{20--25} times as much work as for rank~11.

Since rational elliptic curves are modular
\cite{wiles,taylor-wiles,diamond,cdt,bcdt},
the tables of Cremona \cite{cremona} are complete for~$N\le 20000$.
Hence the lowest conductors for ranks 0--3 are respectively
$11$, $37$, $389$, and $5077$. The rank~4 record was
found by Connell and appears in his Maple\texttrademark\ 
package APECS~\cite{apecs}; the curve has
$[a_1,a_2,a_3,a_4,a_6]=[1,-1,0,-79,289]$
(see Section \ref{model} for notation) and its
conductor of 234446 is more than twice as small as the best example
in~\cite{brumer-mcguinness}. Stein, Jorza, and Balakrishnan
have verified \cite{sjb} that there is no rank 4 elliptic curve of
prime conductor less than 234446.

The rest of this paper is organized as follows.  In the next section
we describe the methods we used to search efficiently for curves
with many small integral points.  We then report on the curves
of low conductor and/or absolute discriminant that we found,
and compare them with previous records.
The next section reports on our computation of further integral points
on each of these record curves and on many others found
in our search.  Finally we compare our numerical results
with previous speculations on the growth of the minimal~$N$\/
as a function of~$r$.

\section{Algorithms}
We describe two algorithms that each find elliptic curves with numerous
integral points whose \hbox{$x$-coordinates} have small absolute value.
The input to our algorithms is an ordered
triple $(h,I,b_2)$ where $h$ is a height parameter,
$I$ is a lower bound on the number of integral points we want,
and $b_2\in\{-4,-3,0,1,4,5\}$,
these being the possible values of $b_2 = a_1^2 + 4 a_2$
for an elliptic curve in minimal Weierstrass form (see below).
We then try to find
elliptic curves $E$ with an equation $y^2=4x^3+b_2x^2+2b_4x+b_6$
such that there are at least $I$ integral points on $E$
with $0\le y\le 2h^3$, $|x|\le h^2$, and $|2b_4|\le 4h^4$.
In modifications of the algorithm, we use variants of these bounds,
and in general only have a high probability of finding
the desired curves.

\subsection{First Algorithm\label{model}}
An elliptic curve $E/\Q$ can be written in its minimal Weierstrass form
as $Y^2+a_1XY+a_3Y=X^3+a_2X^2+a_4X+a_6$, where $a_1$ and $a_3$ are 0~or~1
and $|a_2|\le 1$. We can obtain the ``2-torsion'' equation\footnote{
  The second-named author suggests this term because such a model
  makes it easy to locate the \hbox{$2$-torsion} points on~$E$\/:
  $(x,y) \in E[2]$ if and only if $y=0$.}
$$y^2=4x^3+b_2x^2+2b_4x+b_6$$ by completing the square
via $y=2Y+a_1X+a_3$ and $x=X$, so that we get
$b_2=a_1^2+4a_2$, $b_4=a_1a_3+2a_4$, and $b_6=a_3^2+4a_6$.
Note that this transformation preserves integral points;
we use the 2-torsion equation rather than the minimal equation
since it is relatively fast to check whether its right-hand side is square.
Fixing a choice of $b_2\in\{-4,-3,0,1,4,5\}$ and a height-bound~$h$,
we search for curves with integral points by looping over the
coordinates of such points. In particular, we first fix a $b_4$-value
with $|2b_4|\le 4h^4$ and then loop over integral values of $x$ and $y$
with $|x|\le h^2$ and $0\le y\le 2h^3$, and finally
calculate the value of $b_6$ from the above equation,
counting how often each $b_6$-value occurs.
Note that the above bounds imply that $|y^2|$, $|4x^3|$, and $|2b_4x|$
are all bounded by~$4h^6$.

This algorithm takes on the order of $h^9$ time,
with memory requirements around $h^5$ for the recording of the $b_6$-values.
There are various methods of speeding this up.
We can note that neither positive $b_4$ nor negative~$b_6$
are likely to give curves with many integral points,
due to the shape of the cubic.
From Table~\ref{tbl:relations} we see that $b_4$ cannot
be odd when $b_6$ is even. Also, we know that $b_6$ is a square modulo~4.
We can extend this idea to probabilistic considerations;
for instance, a curve with $b_2=1$ is not that likely to have numerous
integral points unless $b_4$ is odd and $b_6$ is 1~mod~8.
We ran this algorithm for $h=20$, and an analysis showed that
the congruence restrictions most likely  to produce good curves had
$(b_2,b_4,b_6)\,\,{\rm mod}\,\,8$
equal to one of $(1,1,1)$, $(1,3,1)$, $(5,2,4)$, $(5,0,0)$,
$(0,2,1)$, $(0,0,0)$, or $(4,0,1)$.
Of course, there are curves that have many integral
points yet fail such congruence restrictions,
but the percentage of such is rather low (only 10--20\%),
and even those that do have numerous integral points
appear less likely to have high rank. However, our table
of records does contain some curves that fail these congruence
restrictions, so there is some loss in making them.
With these congruence restrictions,
our computation took 15--20 hours on an Athlon~MP~1600
to handle one $b_2$-value for $h=20$;
with no congruence restrictions, this would be about 5~days.
Note that our congruence restrictions imply that the trials
for $b_2=\pm 4$ should only take half as long as the others.
With this algorithm, we broke the low-conductor records of Tom Womack 
(from whose work this sieve search was adapted) for ranks 6, 7, and 8
(see Tables~\ref{tbl:records-cond} and~\ref{tbl:record-compare}).\footnote{
In our tables the stated value of the ``rank'' is actually the rank of
the subgroup generated by small integral points on the curve, which is
very likely to be the actual rank, though in general
such results can be quite difficult to prove.}

\begin{table}[h]
\caption{Congruence relations with $a_1$ and $a_3$
\label{tbl:relations}}
\begin{center}
\begin{tabular}{|c|c|c|c|c|}\hline
$\>a_1\>$&$\>a_3\>$&$\>b_4\>$&$\>b_6\>$&$x$ and $y$\\\hline
$0$&$0$&even&even&$y$ even\\
$0$&$1$&even&odd&$y$ odd\\
$1$&$0$&$\>$even$\>$&$\>$even$\>$&$\>y\equiv x$ (2)$\>$\\
$1$&$1$&odd&odd&$y\not\equiv x$ (2)\\\hline
\end{tabular}
\end{center}
\end{table}

\subsection{Second Algorithm}
The number of elliptic curves~$E$\/ with $b_4 \ll h^4$ and $b_6 \ll h^6$
grows as~$h^{10}$.  The typical such curve has no small integral points
at all: as we have seen, the number of $(E,P)$, with $E$\/ as above
and $P \in E(\Q)$ a small integral point, grows only as $h^9$,
as does the time it takes to find all these $(E,P)$.
But we expect that even in this smaller set the typical~$E$\/
does not interest us, because it has no integral points
other than $\pm P$.
We shall see that there are (up to at most a logarithmic factor)
only $O(h^8)$ curves~$E$\/ in this range together with a {\bf pair}
of integral points $P,P'$ such that $P' \neq \pm P$,
and that again we can find all such $(E,P,P')$ with given $b_2,b_4$
in essentially constant time per curve.
We thus gain a factor of almost~$h$\/ compared to our first
algorithm.\footnote{This pair-finding idea is also used
in~\cite{elkies-rogers} to find curves $x^3+y^3=k$ of high rank.}
Further improvements might be available by searching
for elliptic curves with three or more points, but we do not know
how to do this with the same time and space efficiency.

We wish to compute all $b_4,b_6,x_1,y_1,x_2,y_2$ in given ranges
that satisfy the pair of equations
$y_j^2 = 4 x_j^3 + b_2 x_j^2 + 2 b_4 x_j + b_6$ $(j=1,2)$.
Subtracting these two equations, we find that
$$(y_2-y_1)(y_2+y_1)=(x_2-x_1)[2b_4+b_2(x_2+x_1)+4(x_2^2+x_1x_2+x_1^2)].$$
We can thus write
$$x_2-x_1=rt\qquad{\rm and}\quad y_2-y_1=rs\qquad{\rm and}\quad y_2+y_1=tu$$
for some integers $r,s,t,u$. From the latter two equations,
we see that we need $rs$~and~$tu$ to have the same parity
in order for the $y$'s to be integral.
Our expectation is that generically we shall have
$r,t\ll h$ and $s,u\ll h^2$ when the $x$-values are
bounded by $h^2$ and the $y$-values by~$h^3$.
It is unclear how often this expectation is met.
One way of estimating the proportion is to consider pairs
of points $(x_1,y_1),(x_2,y_2)$ with $|x_i|\le h^2$ on various curves
and see what values of $(r,s,t,u)$ are obtained. This is not quite
well-defined from the above; for instance, the quadruple
$(x_1,x_2,y_1,y_2)=(7,3,6,2)$ could have $(r,s,t,u)$ as either
$(4,1,1,8)$ or $(2,2,2,4)$. However, it becomes well-defined
upon imposing the additional condition that $r=\gcd(y_2-y_1,x_2-x_1)$.
Experiments show that about 18\% of the $(r,s,t,u)$ obtained from
this process satisfy $1\le r,t\le h$, though the exact percentage can
vary significantly with the curve.
Note that swapping $r$ and $t$ or negating either
leads either to a switching of $(x_1,y_1)$ and $(x_2,y_2)$
or to a negation of $y$-values. Thus we can assume that $1\le r\le t$.

We rewrite the above equation in the form
$$rstu=rt[2b_4+b_2(x_1+x_2)+3(x_1+x_2)^2+(x_1-x_2)^2]$$
and define $z=x_1+x_2$ so that $su=2b_4+b_2z+3z^2+(rt)^2$.
Our algorithm is now the following.
Given one of the six possible values of~$b_2$,
we loop over $2b_4$-values between $-4h^4$ and $0$
(implementing our above comment that positive $b_4$-values
are not that likely to give curves with many integral points).
For each value of $b_4$ we loop over pairs of integers $(r,t)$
that satisfy $1\le r\le t\le h$.
We then compute $l=rt$ and loop over values of $2x_2$ (that is, $z+l$)
with $-2h^2\le 2x_2\le 2h^2$.
Next we compute the quantity $W=2b_4+b_2z+(l^2+3z^2)$ and factor this in
all possible ways as $W=su$. We then take $y_2=(rs+tu)/2$
(assuming that $rs$ and $tu$ have the same parity)
and compute $b_6=y_2^2-4x_2^3-b_2x_2^2-2b_4x_2$.
As before, we record the $b_6$-values and count how many times each occurs.
This algorithm takes about $h^8\log h$ time, where the logarithmic
factor comes from solutions of $W=su$, assuming we can find these
relatively fast via a lookup table. Already at $h=20$, a version
of this algorithm ran in under an hour and found most of the curves
found by the first algorithm. One can view this algorithm
as looping over pairs of $x$-values (both of size $h^2$),
or more precisely the sum (given by~$z$) and difference (given by~$l$)
of such a pair, and then reconstructing the $y$-values by factoring.
Thus the inner loop takes time $h^4\log h$
instead of $h^5$ as in the first algorithm.

\subsection{Implementation Tricks}
We now describe the various tricks we used in the implementation.
We shall see that our $b$-congruence restrictions allow us to limit
the $z$ and $l$ values in a productive way.
First we consider the cases where $b_2$ is odd.
Given a fixed $b_4$-value we only loop over $z$'s and $l$'s that are both odd,
and can note that this makes $W$ odd.
Actually we do not loop over~$l$ but determine it as $l=rt$;
thus we are looping over odd $r$ and $t$ with $1\le r\le t\le h$.
It may seem that this loses a factor of 4 of $(r,s,t,u)$ quadruples
(with the $z$-restriction losing nothing because $z$ must be of the same
parity as~$rt$), but we claim that it is actually only a
factor of 2 for ``interesting'' curves. Indeed, though our yield of
$b_6$-values will drop by a factor of 4 because of this parity
restriction on both $r$ and~$t$, many of these values of $b_6$
will correspond to curves on which all integral points
have $x$-coordinates of the same parity.  Since $l$ is the difference of two
$x$-coordinates, this implies that $l$ must be even for all pairs of
integral points. These curves, which are plainly less likely to have
a large number of integral points, are over-represented in the curves
we ignore through not considering even~$l$. From this we get our
heuristic assertion that restricting to odd~$l$ loses only
a factor of about~2.

We also consider only the values of~$z$
for which $|W|$ is less than a certain bound.
This serves a dual purpose in that it speeds up the algorithm
and also reduces the size of the tables used for factoring.
We see that $W=2b_4+b_2z+(l^2+3z^2)$
should be of size~$h^4$, and so we restrict the size of $W$ via the
inequality $|W|\le 2h^4/U$, where $U$ is a parameter we can vary
(we had $U=1$ for the experiments with $h=30$ and~$h=40$).
Again it is not immediately clear how many $(r,s,t,u)$ quadruples we
miss by making this restriction on~$W$, and again the proportion
can depend significantly upon the curve (curves with $b_4$ near $-2h^4$
lead to more quadruples with large $|W|$ than those with $b_4$ close to~0).
Experimentation showed that with $U=1$ we catch on average about
$83\%$ of the relevant $b_6$-values under this restriction.
Our expectation might be approximate inverse linearity of the catch rate
in~$U$, though only in the limit as $U\rightarrow\infty$.
Experimentation showed that with $U=8$ our catch rate is down to 27\%,
while at $U=32$ it is about 10\%.
However, there is interdependence between this restriction and that
on the size of $r$ and $t$ --- when $r$ and $t$ are both small, this
corresponds to a small \hbox{$x$-difference,}
which implies a small \hbox{$y$-difference,}
and so $W=su$ should also be diminished in size.
In a final accounting of the proportion of $(r,s,t,u)$ quadruples,
including the loss of a factor of~2 from the parity restriction on $l$,
we find that with $U=1$ we catch 7.4\% of the quadruples,
with $U=8$ we catch 2.5\%, and with $U=32$ we catch just under 1\%.
Most of the curves of interest to us
have at least 40 integral points within the given bounds $|x|\le h^2$
and $0\le y\le 2h^3$, and thus have at least 780 pairs of integral points.
So missing 93\% or more of the $(r,s,t,u)$ quadruples does not
trouble us --- indeed, our ``laziness'' in not finding all the
possible $(r,s,t,u)$ is quadratically efficient compared to what we
would achieve via similar ``laziness'' in our first algorithm.

So far we assumed that $b_2$ was odd, but similar ideas
apply also for even values of~$b_2$.
When $b_2=\pm4$, we took $l$ and $z$ to be even but not congruent modulo~4.
This ensures that $W$ is 4~mod~8.
Similarly, when $b_2=0$ and $b_6$ is odd,
we take $l$ and $z$ to be even and congruent modulo~4,
again ensuring that $W$ is 4~mod~8.
To implement these restrictions, we took $r \equiv 2 \bmod 4$
with no restriction on $t$ other than $t \geq r$ if $t$ itself is also 2~mod~4
(with both variables less than $h$ as before). Again we required
that $|W|\le 2h^4/U$, and here we have various restrictions on
the decomposition $W=su$ depending on $t$~mod~4.
Specifically, we can always take $s$ odd and $4|u$,
we can take both $s$ and $u$ even if $t$ is odd,
and we can take $4|s$ and $u$ odd if $t$ is 2~mod~4,
as we need for $y$ to be odd in these cases.
When $b_2=0$ and $b_6$ is even, we take $z$ and $l$ both to be odd,
which makes $W$ be 4~mod~8, and we need both $s$ and $u$ to be even
(hence each 2~mod~4) for $y$ to be even. As above, the loss in
the number of interesting $b_6$-values from these restrictions
is not much more than a factor of~2.

\subsection{More Tricks}
To reduce the memory needed for the counting of $b_6$ values,
we used the following idea.
We create a array of $2^L$ counters (of size $16$ bits each);
for instance, for $h=30$ we used $L=19$.
Then for each $b_6$-value we obtain from above,
we reduce $\lfloor b_6/8\rfloor$ modulo $2^L$,
and increment the corresponding counter.
In other words, we only record $b_6$ modulo $2^{L+3}$.
At the end of the loops over $r$, $t$, and $z$,
we extract the counters with at least 10 hits.
These residue classes are then passed to a secondary test phase.
Here we set up counters for the values of $b_6$ with
$0\le b_6\le 4h^6$ that are in
the desired residue class $b_6^\#$ modulo $2^{L+3}$.
We then run through integral $x$ with $|x|\le h^2$,
and for each \hbox{$x$-value} determine
the corresponding positive $y$-values such that
$y^2\equiv 4x^3+b_2x^2+2b_4x+b_6^\# \bmod{2^{L+3}}$
via a lookup table of square roots modulo $2^{L+3}$.
Most of these \hbox{$y$-values} exceed~$2h^3$,
and we thus ignore them. If not, we compute $b_6$ exactly
from $x$ and $y$, and increment the corresponding counter.
After running over all the \hbox{$x$-values,}
we then check for large counter values.
By taking $2^L$ somewhere around~$h^4$
(note that this is about how many \hbox{$b_6$-values} we generate),
we can use this method to handle a \hbox{$b_6$-congruence-class}
in essentially $h^2$ time.
This is generically small compared to the $h^4$ time for the loops
over $r$, $t$, and $z$; when $h=30$ we averaged about 100 congruence
classes checked for each \hbox{$b_4$-value,}
but the time for the loops still dominated.

\section{Experimental Results}
We ran this algorithm with $h=30$ and $U=1$
with a few more congruence classes in consideration,
taking about a day for each $(b_2,b_4,b_6)$ class.
We then proceeded to run it for $h=40$ and $U=1$,
and then $h=60$ and $U=8$, taking a few weeks for each $(b_2,b_4,b_6)$ class.
Other runs were done with the ``better'' congruence restrictions
of $(b_2,b_4,b_6)$ with varied parameters up to $h=90$ and $U=48$.
Though with the $U=48$ restriction we are catching less than 1\%
of the $(r,s,t,u)$ quadruples, by this time we expect that most interesting
curves have 60 or more integral points with $|x|\le h^2$ and $0\le y\le 2h^3$;
indeed, even with the $h=60$ search all the record curves
we found had at least 70 integral points in this range.

\begin{table}
\caption{Low conductor records for ranks 5--11\label{tbl:records-cond}}
\begin{center}
\begin{tabular}{|@{\,}l@{\,}|@{\,}r@{\,}|r|c|c|}\hline
$[a_1,a_2,a_3,a_4,a_6]$ & $N$ & $|\Delta|/N$ & $I$ & $r$\\\hline
$[0,0,1,-79,342]$&19047851&1&39&5\\
$[1,0,0,-22,219]$&20384311&1&29&5\\
$[0,0,1,-247,1476]$&22966597&1&40&5\\
$[1,-1,0,-415,3481]$&34672310&10&52&5\\
$[0,0,0,-532,4420]$&37396136&32&52&5\\\hline

$[1,1,0,-2582,48720]$&5187563742&6&71&6\\
$[0,0,1,-7077,235516]$&5258110041&243&67&6\\
$[1,-1,0,-2326,43456]$&5739520802&2&60&6\\
$[1,-1,0,-16249,799549]$&6601024978&184&68&6\\
$[1,-1,1,-63147,6081915]$&6663562874&32768&88&6\\\hline

$[0,0,0,-10012,346900]$&382623908456&32&101&7\\
$[1,0,1,-14733,694232]$&536670340706&8&77&7\\
$[0,0,1,-36673,2704878]$&814434447535&5&84&7\\
$[1,-1,0,-92656,10865908]$&858426129202&142&92&7\\
$[1,-1,0,-18664,958204]$&896913586322&26&109&7\\\hline

$[1,-1,0,-106384,13075804]$&249649566346838&14&124&8\\
$[1,-1,0,-222751,40537273]$&292246301470558&2&101&8\\
$[0,0,0,-481663,128212738]$&314214346667560&160&141&8\\
$[1,-1,0,-71899,5522449]$&314658846776578&34&130&8\\
$[1,-1,0,-124294,14418784]$&315734078239402&106&131&8\\\hline

$[1,-1,0,-135004,97151644]$& 32107342006814614 &122&191&9\\
$[1,-1,0,-613069,98885089]$& 43537345103385386 &242&203&9\\
$[0,0,1,-3835819,2889890730]$& 62986816173592807 &67&142&9\\
$[1,0,1,-1493028,701820182]$& 72070075910145406 &2&139&9\\
$[1,0,1,-1076185,496031340]$& 77211251506212554 &344&156&9\\\hline

$[0,0,1,-16312387,25970162646]$& 10189285026863130793&1331&262&10\\
$[1,-1,0,-10194109,12647638369]$& 22006161865320788846& 58&241&10\\
$[0,0,1,-21078967,35688990786]$& 22630148490190627609 &2173 &238&10\\
$[1,-1,0,-1536664,648294124]$& 25440555737235843986 &2&207&10\\
$[1,-1,0,-4513546,3716615296]$& 39432942782223365758 &2&179&10\\\hline

$[0,0,1,-16359067,26274178986]$& 18031737725935636520843& 1&229&11\\
$[1,-1,0,-38099014,115877816224]$& 66484354768372183177742& 34&281&11\\
$[1,-1,0,-41032399,106082399089]$& 219576020293485812169274& 2&236&11\\
$[1,-1,0,-34125664,69523358164]$& 227946110025657660240686& 2&215&11\\
$[1,-1,0,-56880994,168642718624]$& 252948166615918192888894& 2&235&11\\\hline
\end{tabular}
\end{center}
\end{table}

\begin{table}
\caption{Value of $\log N$ for
old and new rank records\label{tbl:record-compare}}
\begin{center}
\begin{tabular*}{3.0in}{@{\extracolsep{\fill}}|c|c|c|c|c|c|c|}\hline
&6&7&8&9&10&11\\\hline
Old&22.383&27.703&33.962&40.721&49.033&55.852\\
New&22.370&26.670&33.151&38.008&43.768&51.246\\\hline
\end{tabular*}
\end{center}
\end{table}

Table \ref{tbl:records-cond} lists minimal equations for each of the
five curves of smallest conductor $N$ for each rank from 5--11
that were found by the above method.
Table \ref{tbl:records-disc} lists similar data for
smallest absolute discriminant~$|\Delta|$.
The rank 5 data agree with the data from the
Elliptic Curve Database~\cite{stein-watkins}.
The $I$ column gives how many $x$-coordinates of integral points
we found (see Section~\ref{intpts}) for the given equation.
Some of the curves
fail our congruence conditions on~$(b_2,b_4,b_6)$,
but they still can be found via a non-minimal model;
indeed, letting $c_4$ and $c_6$ be the invariants of the minimal model,
the model with invariants $12^4c_4$ and $12^6c_6$
has $b_2=0$ and $4|b_4$ and $8|b_6$ and is thus in the $(0,0,0)$ class.
In this way, from $(b_2,b_4,b_6)=(0,-1826496,2637633024)$
we recover the curve $[1,0,0,-22,219]$.

\begin{table}[h]
\caption{Low absolute discriminant
records for ranks 5--10\label{tbl:records-disc}}
\begin{center}
\begin{tabular}{|@{\,}l|r@{\,}|@{\,}c@{\,}|c|}\hline
$[a_1,a_2,a_3,a_4,a_6]$ & $|\Delta|$ & $I$ & $r$\\\hline
$[0,0,1,-79,342]$&19047851&39&5\\
$[1,0,0,-22,219]$&20384311&29&5\\
$[0,0,1,-247,1476]$&22966597&40&5\\
$[0,1,1,-100,110]$& 55726757&33&5\\
$[0,0,1,-139,732]$& 59754491&32&5\\\hline

$[1,0,0,-9227,340354]$& 6822208199&36&6\\
$[0,0,1,-277,4566]$& 7647224363&49&6\\
$[0,0,1,-379,5172]$& 8072781371&51&6\\
$[0,0,1,-889,9150]$& 8796007189&54&6\\
$[0,1,1,-390,5460]$& 9694585723&43&6\\\hline

$[0,0,1,-1387,68046]$& 1829517077483&71&7\\
$[0,0,1,-5707,151416]$& 1991659717477&68&7\\
$[1,0,1,-5983,164022]$& 2010552189452&72&7\\
$[1,0,1,-14505,667472]$& 2132568452204&71&7\\
$[0,0,1,-15577,744876]$& 2206378706437&71&7\\\hline

$[0,1,1,-23846,1022562]$& 409086620841461&78&8\\
$[0,0,1,-23737,960366]$& 457532830151317&96&8\\
$[0,1,1,-16440,1394010]$& 561715239383323&84&8\\
$[1,-1,0,-222751,40537273]$& 584492602941116&101&8\\
$[1,-1,0,-201814,34925104]$& 643509175703572&109&8\\\hline

$[0,0,1,-167419,30261330]$& 95276302704064331&135&9\\
$[1,0,1,-1493028,701820182]$& 144140151820290812&139&9\\
$[0,0,1,-514507,140806716]$& 151673348057775877&126&9\\
$[0,0,1,-402157,96291336]$& 157107745029925477&131&9\\
$[0,0,1,-826609,289956150]$& 172539371946838571&120&9\\\hline

$[1,-1,0,-1536664,648294124]$&50881111474471687972&207&10\\
$[0,0,1,-1788817,843180666]$& 59202439687694448757&176&10\\
$[1,-1,0,-4513546,3716615296]$& 78865885564446731516&179&10\\
$[0,1,1,-1856500,1072474760]$& 87950374485438204043&154&10\\
$[0,0,1,-2438527,1545098346]$& 103294665688000244363&173&10\\\hline
\end{tabular}
\end{center}
\end{table}

How good is this method at finding elliptic curves of low
conductor~$N$\/ and relatively high rank?  Obviously if such a curve
has few integral points then we will not find it.
Indeed, it was suggested to us by J.~Silverman that for large ranks~$r$
the smallest conductor curve might not have~$r$
independent integral points.
However, for the ranks we consider
there are sufficiently many independent integral points;
the same is true for Mestre's rank 15 curve~\cite{mestre15},
but apparently not for later rank records.

Note also that our search
operates by increasing $b_4$ and $b_6$ corresponding to some height
parameter, which is not exactly the same as simply increasing the
absolute value~$|\Delta|$ of the discriminant,
which again is not quite the same as just increasing~$N$.
Finally, the probabilistic nature of our
algorithm and the necessity of restricting to ``likely'' congruence
classes also cast doubt on the exhaustiveness of our search procedure.
However, we are still fairly certain that the curves we found for ranks 5--8
are indeed the actual smallest conductor curves for those ranks.
Note that our methods were almost exhaustive in the region of interest
($h$~up to about 30), and were verified with the first algorithm in
much of this range. For rank 9 we could be missing some curves with
large $|\Delta|/N$ and $h$ around 50 or so; perhaps this range should
be rechecked with a smaller $U$-parameter. Indeed,
for a long time the second curve on the $r=10$ list,
which we found with an $h=60$ search,
was our rank~$10$ record, but then a run with $h=80$ found the first
and third curves which have $h$-values of about 64 and 68 respectively.
We have yet to find many rank~11 curves with large $|\Delta|/N$;
note that our current record curve in fact has prime conductor.
This suggests that there still could be significant gains here.
However, Table \ref{tbl:record-compare}, which lists values of $\log N$
for the old records and our new ones, indicates that our method has
already shown its usefulness.
There does not seem to have been any public compilation of such records before
Womack~\cite{womack} did so on his website in the year 2000, soon after
he had found the records for ranks 6--8 via a sieve-search.\footnote{
Prior to Womack, there were records listed in the Edinburgh dissertation
of Nigel Suess (2000); it appears that Womack and Suess enumerated these
lists in part to help dispense Cremona from answering emails about the
records. Indeed, the rank~4 record of McConnell \cite{apecs} mentioned above
was relatively unknown for quite a while.}
The records for ranks 9--10 were again due to Womack but from a
Mestre-style construction~\cite{mestre-method-1991},
with Mestre listing the rank~11 record in \cite{mestre}
(presumably found by the methods of~\cite{mestre-method-1982}).
There was no particular reason to
expect the old records for ranks 9--11 to be anywhere near the true
minima, as they were constructed without a concentrated attempt
to make the conductor or any related quantity as small as possible.

\section{Counting Integral Points}\label{intpts}
\vspace*{-1.5pt}
Since we have curves that have many integral points of small height,
it is natural to ask how many integral points these curves have overall,
with no size restriction.
For our curves of rank higher than~8,
current methods, as described in~\cite{ST},
do not yet make it routine to list
all the integral points and to prove that the list is complete.
Indeed, even verification that the rank is actually what
it seems to be is not necessarily routine.

However, we have at least two ways to search for integral points
and thus obtain at least a lower bound for the number of integral points.
One method is a simple sieve-assisted search, which can reach
$x$-values up to~$10^{12}$ in just over an hour on an Athlon~MP~1600.
The other method is to write down a linearly independent set from the
points we have, and then take small linear combinations of these.\footnote{
The possible size of coefficients in such linear combinations can be bounded
via elliptic logarithms (possibly $p$-adic) as in \cite{ST} and later works.
Also, as indicated by Zagier~\cite{zagier}, one can combine elliptic logarithms
with lattice reduction to search for large integral points, but we did
not do this.}
For this, we took the linearly independent set
that maximized the minimal eigenvalue of the height-pairing matrix
(as in the ``$c_1$-optimal basis'' of~\cite{ST}),
subject to the condition that the set must
generate (as a subgroup of $E(\Q)$) all the integral points in our list.
We then tried all $((2m+1)^r-1)/2$ relevant linear combinations with
coefficients bounded in absolute value by~$m$.\footnote{
We do not compute sums of points on elliptic curves directly over the
rationals, but instead work modulo a few small primes
and use the Chinese Remainder Theorem.}
With $r=11$ and $m=3$ this takes about an hour.

The maximal number of integral points
we found on a curve was $281\times 2$ for~$[1,-1,0,-38099014,115877816224]$.
This can be compared with the rank 14 curve
$[0,0,1,-2248232106757,1329472091379662406]$ that is
listed by Mestre \cite{mestre},
which we find has $311\times 2$ integral points with $|x|\le 10^{12}$,
plus at least $7\times 2$ more that were found with linear combinations
as above. Note that amongst the curves of a given rank there is not much
correspondence between number of integral points and smallness
of conductor. For instance, we have no idea which curve of rank~7
has the maximal number of integral points;\footnote
{The maximal count of integral points may well be attained by a curve
with nontrivial torsion, whose discriminant would then be too large
to be found by our search.}
trying a few curves
found in our search turned up $[1,-1,0,-22221159,40791791609]$
which has at least $165\times 2$ integral points, but conductor
$13077343449126$, more than $34$ times larger than our record.
Note that $|\Delta|/N = 2 \cdot 3^4 11^4 23^2$ for this curve;
in general, large values of $|\Delta|/N$\/ seem to correlate
with large counts of integral points (see Table \ref{tbl:records-cond}).

\section{Growth of Maximal Rank as a Function of Conductor}\label{rankcond}
\vspace*{-0.5pt}
We review two different heuristics and conjectures for the growth of the
maximal rank of an elliptic curve as a function of its conductor,
and then indicate which is more likely according to our data.
The first conjecture is due to Murty and appears in the
appendix to~\cite{rajan}. He first notes that,
similar to a heuristic of Montgomery \cite[pp.~512--513]{montgomery}
regarding the \hbox{$\zeta$-function,} it is plausible that
$\arg L_E(1+it)\ll\sqrt{\log(Nt)/\log\log(Nt)}$ as $t\rightarrow\infty$.
Murty speculatively applies this bound in a small circle of radius
$1/\log\log N$ about $s=1$.
He then claims that Jensen's Theorem implies that
the order of vanishing of $L_E(s)$ at $s=1$ is bounded
by $C\sqrt{\log N/\log\log N}$,
though we cannot follow the argument.
Assuming the Birch--Swinnerton-Dyer conjecture~\cite{bsd},
the same upper bound holds for the rank of the elliptic curve.
However, the Montgomery heuristic comes from taking
the approximation $\log\zeta(s)=\sum_{p\le t} p^{-\sigma-it}+O_\sigma(1)$
for~$\sigma>1/2$
and assuming that the $p^{-it}$ act like random variables; upon taking
a limit as $\sigma\rightarrow 1/2$, this implies the asserted bound
of $\sqrt{\log t/\log\log t}$, but only for large~$t$.
Indeed, in our elliptic curve case with small~$t$,
we should have an approximation (see \cite{goldfeld})
more like $\log L_E(s)\sim\sum_{p\le X} a_p/p^{-s}$;
it is unclear whether the variation of the $a_p$ or that of $p^{-s}$
should have the greater impact. Finally, Conrey and Gonek \cite{conrey-gonek}
contest that Montgomery's heuristic could be misleading;
they suggest that $\log |\zeta(1/2+it)|$
(and maybe analogously the argument)
can be as big as $C\log t/\log\log t$ instead of the square root of this.
One idea is that the above limit as $\sigma\rightarrow 1/2$
disregards a possibly larger error term coming from zeros of~$\zeta(s)$;
the asymmetry of upper and lower bounds
for $\log |\zeta(1/2+it)|$ makes the analysis delicate.

A classic paper of Shafarevich and Tate \cite{sha-tate} shows
that in function fields the rank grows at least as fast as 
the analogue of ${\log N\over 2\log\log N}$.
However, the curves used in this construction were isotrivial, and thus
fairly suspect for evidence toward a conjecture over number fields.
Ulmer recently gave non-isotrivial function field examples
with this growth rate, and conjectured \cite[Conjecture~10.5]{ulmer}
that this should be the proper rate of growth even in the number field case,
albeit possibly with a different constant.
In a different paper \cite[p.~19]{ulmer2}, Ulmer notes that
certain random matrix models suggest that the growth rate is as
in the function field case; presumably this is an elliptic
curve analogue of the work of~\cite{conrey-gonek}.

\begin{figure}
\begin{center}
\scalebox{1.0}{\includegraphics{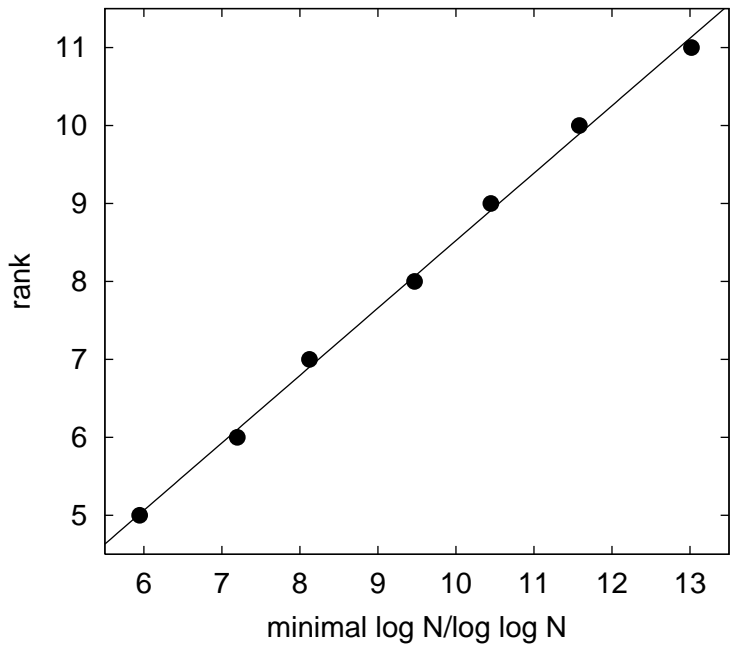}}
\caption{Plot of rank versus $\log N/\log\log N$\label{fig:graph}}
\end{center}
\end{figure}

Figure \ref{fig:graph} plots the rank $r$ versus $\log N/\log\log N$,
where $N$ is the smallest known conductor for an elliptic curve of rank~$r$.
A log-regression gives us that the best-fit exponent is~0.975,
much closer to the exponent of 1.0 of Ulmer than to the 0.5 of Murty.
Note that an improvement in the records for ranks 9--11 would most
likely increase the best-fit exponent.

\eject

\vspace*{-4ex}

Assuming the growth is linear, the line of best fit is approximately
$r=0.865{\log N\over\log\log N}-0.126$.
But this could mislead; GRH plus BSD implies
$$r\le {1\over 2}{\log N\over\log\log N}
\biggl(1+{\log 8e\over\log\log N}
+O\biggr({1\over(\log\log N)^2}\biggr)\biggr).$$
The main term in the above already appears in Corollary~2.11 of~\cite{brumer}
(see also Proposition~6.11),
and we have simply calculated the next term in the expansion.
To get more reliable data, we would need to consider curves
with $\log\log N$ rather large, which is of course quite difficult.

Finally we mention a possible heuristic refinement of the above upper bound.
The bound comes from a use of the Weil explicit formula
(see \cite[2.11]{brumer}) to obtain the relation
$\sum_\gamma h(\gamma\log\log N)\sim {\log N\over 2\log\log N}$,
where $h(t)=\bigl({\sin t\over t}\bigr)^2$ and
the sum is over imaginary parts of nontrivial zeros of~$L_E(s)$,
counted with multiplicity.
When only the high-order zero at $\gamma=0$ contributes,
we get the stated upper bound. In the function field case,
the other zeros contribute little because they are all near the minima of~$h$.
This is unlikely to occur in the number field case.
Also unlikely is the idea that the other zeros have negligible
contribution due to the $1/t^2$ decay of~$h$.
Thus the other zeros are likely to have some impact;
however, it is not clear how large this impact will be.

\end{document}